\documentstyle {article}

\title{\bf COGROWTH AND ESSENTIALITY IN GROUPS AND ALGEBRAS}

\author
{Amnon Rosenmann\thanks{Supported by the Minerva Fellowship} \\
\vspace {-2 mm}
\small
Institute for Experimental Mathematics \\
\vspace {-2 mm}
\small
University GH Essen \\
\vspace {-2 mm}
\small
45326 Essen, Germany \\
\vspace {-2 mm}
\small
amnon@exp-math.uni-essen.de
}

\date{}

\begin{document}
\maketitle

\newtheorem{clm}{Claim}[section]
\newtheorem{rmrk}[clm]{Remark}
\newtheorem{exm}[clm]{Example}
\newtheorem{dfn}[clm]{Definition}
\newtheorem{thm}[clm]{Theorem}
\newtheorem{lem}[clm]{Lemma}
\newtheorem{prop}[clm]{Proposition}
\newtheorem{corol}[clm]{Corollary}

\newcommand{\FF} {\mbox{$\cal F$}}
\newcommand{\ds} {\mbox{$\scriptscriptstyle S$}}
\newcommand{\dt} {\mbox{$\scriptscriptstyle T$}}
\newcommand{\dghAB} {\mbox{$\scriptscriptstyle G / H_{1} \cap H_{2}$}}
\newcommand{\dghB} {\mbox{$\scriptscriptstyle G/H_2$}}
\newcommand{\dghA} {\mbox{$\scriptscriptstyle G/H_1$}}
\newcommand{\dgh} {\mbox{$\scriptscriptstyle G/H$}}
\newcommand{\dg} {\mbox{$\scriptscriptstyle G$}}
\newcommand{\dh} {\mbox{$\scriptscriptstyle H$}}
\newcommand{\driAB} {\mbox{$\scriptscriptstyle R/I_1 \cap I_2$}}
\newcommand{\driB} {\mbox{$\scriptscriptstyle R/I_2$}}
\newcommand{\driA} {\mbox{$\scriptscriptstyle R/I_1$}}
\newcommand{\drii} {\mbox{$\scriptscriptstyle rR+I/I$}}
\newcommand{\dri} {\mbox{$\scriptscriptstyle R/I$}}
\newcommand{\dr} {\mbox{$\scriptscriptstyle R$}}
\newcommand{\di} {\mbox{$\scriptscriptstyle I$}}
\newcommand{\lrw}{\mbox{$\longrightarrow$}}
\newcommand{\Lrw}{\mbox{$\Longrightarrow$}}
\newcommand{\Llrw}{\mbox{$\Longleftrightarrow$}}
\newcommand{\Notin}{\mbox{$\in \! \! \! \! \! / \,$}}

\begin{abstract}
The cogrowth of a subgroup is defined as the growth of a
set of coset representatives which are of minimal length.
A subgroup is essential if it intersects non-trivially every
non-trivial subgroup. The main result of this paper is that
every function $f : \mbox{\bf N} \cup \{ 0 \} \lrw \mbox{\bf N}$
which is strictly increasing, but at most exponential, is
equivalent to a cogrowth function of an
essential subgroup of infinite
index of the free group of rank two. This class of functions
properly contains the class of growth functions of groups.

The notions of growth and cogrowth of right ideals in algebras
are introduced. We show that when the algebra is without 
zero divisors then every right ideal, whose cogrowth is less
than that of the algebra, is essential.
\end{abstract}

\section{Growth, Cogrowth and Essentiality in Groups}
\subsection{Growth and Cogrowth of Subgroups}
A {\em growth function} $\Gamma_{\ds}(n)$ 
on a set $S$ with a length function $l$ on it is defined by
\begin{equation}
\Gamma_{\ds}(n) := \mbox{card} \{ s \in S \ \mid \ l(s) \leq n \},
\end{equation}
assuming that $\Gamma_{\ds}(n)$ is finite for each $n$.
A preorder is given on the growth functions by
\begin{equation}
\Gamma_{1}(n) \preceq \Gamma_{2}(n) \ \Llrw \
\exists C \ [ \Gamma_{1}(n) \leq \Gamma_{2}(Cn) \ ].
\end{equation}
The notion of growth when applied to finitely generated
groups (see \cite{Gri1} for an overview) has been investigated
mainly after Milnor's paper (\cite{Mil}). A geometric 
interpretation can be given, for example, when computing the
growth function of the fundamental group of a Riemannian manifold.
In order to avoid the dependence of the length function upon
the generating set of the group, an equivalence relation 
is used
\begin{equation}
\Gamma_1(n) \sim \Gamma_2(n) \ \Llrw \ \exists C \ [ \Gamma_1(n)
\leq \Gamma_2(Cn) \ \& \ \Gamma_2(n) \leq \Gamma_1(Cn) \ ].
\end{equation}
We will also use the following notation
\begin{equation}
\Gamma_1(n) \prec \Gamma_2(n) \ \Llrw \ \Gamma_1(n) \preceq 
\Gamma_2(n) \ \& \ \Gamma_1(n) \not \sim \Gamma_2(n).
\end{equation}
The growth function of the group $G$ will be denoted by
$\Gamma_{\dg}(n)$, when referring to its equivalence class and also
when some fixed generating set is assumed (but omitted in the
notation).

When $H$ is a subgroup of $G$ we may speak of the {\em cogrowth}
of $H$ in $G$, denoted $\Gamma_{\dgh}(n)$. This is defined to be
the growth of a (complete) set of coset representatives 
for $H$ in $G$ which is ``minimal'' in the sense that every
representative is of minimal length in its coset relative to 
the given group-generating set of $G$. Clearly, any other set of coset
representatives will grow at most as fast as a minimal set. 
We also notice that cogrowth functions of $H$ relative to
different generating sets of $G$ are equivalent. 
(Remark: it does not matter if we take
right or left cosets because the inverses of the right representatives
can be used as left representatives). Constructing the set of
right coset representatives by induction on length, we see that 
it can always be chosen so that its elements are initially
closed, that is arranged in a form of a tree. Such a set is
called a Schreier transversal (see \cite{Lyn}).
Presenting an order on the generating set of $G$
induces a ``ShortLex'' 
total order on $G$ (comparing elements first by length and then
by the lexicographic order). Relative to this order each subgroup
has a unique minimal transversal, which is also a Schreier
transversal.

The growth of a subgroup $H < G$ with respect to the generators
of $G$ was studied by Grigorchuk (see \cite{Gri2}). We denote
this growth function by $\Gamma^{({\dg})}_{\dh}(n)$ (to
distinguish it from the growth function of $H$ when
considered a {\em group}).
A connection between the different growth functions can be given by
(assuming a fixed generating set)
\begin{equation}
\sum_{i=0}^n \gamma_{\dgh}(i) \Gamma^{({\dg})}_{\dh}(n-i) \leq
\Gamma_{\dg}(n) \leq
\sum_{i=0}^n \gamma_{\dgh}(i) \Gamma^{({\dg})}_{\dh}(n+i),
\end{equation}
where $\gamma_{\dgh}(i) := \Gamma_{\dgh}(i) - \Gamma_{\dgh}(i-1)$,
that is the coset representatives of length exactly $i$.
The left inequality comes from the fact that different cosets
are disjoint subsets. The right inequality is
by the definition of the cogrowth function through the coset
representatives of minimal length.

When $H$ is a normal subgroup of $G$ then the cogrowth of $H$ 
describes the
growth of the group $G/H$. This imposes restrictions on the
cogrowth of $H$:
\begin{equation}
\Gamma_{\dgh}(n_1+n_2) \leq \Gamma_{\dgh}(n_1 -1)
+ \gamma_{\dgh}(n_1) \Gamma_{\dgh}(n_2).
\label{eqGroups10}
\end{equation}
Since each right coset is also a left one when $H$ is normal,
the minimal Schreier transversal tree $T$, relative to a
ShortLex order, is suffix as well as prefix-closed. In other words,
each subtree $T'$ of $T$ is ``covered''
by $T$ when putting the root of $T$ over the root of $T'$.
(This is also the reason why a finitely generated subgroup
of a free group which contains a non-trivial normal subgroup 
is of finite index.)
Thus if $g=g_1g_2 \in T$ when written in reduced form then
the subtree $T(g)$ with root in $g$ is ``contained'' 
(in the above sense) in the subtree
$T(g_2)$. Therefore, there exists a descending chain of subtrees of
$T$, which is of length $l(g)$ and
terminates with $T(g)$. This tendency of $T$ to
``close'' itself raised the question if there are groups of
non-exponential but also non-polynomial growth. Milnor
and Wolf showed that a f.g. solvable group has polynomial
growth if it is virtually-nilpotent and otherwise has exponential
growth (\cite{Mil1}, \cite{Wolf}). Gromov (\cite{Gro}) 
showed that virtually-nilpotent groups are the only ones with
polynomial growth. Then
Grigorchuk succeeded to obtain remarkable examples of
groups of ``intermediate growth'' (between $n^d$ and $e^n$), 
and to show that the set 
of growth degrees of finitely generated
groups is of the continuum cardinality
(see \cite{Gri}, \cite{Gri1} and also \cite{Fab}).

If $H_1, H_2$ are subgroups of $G$ then
$\Gamma_{\dghAB}(n) \geq max \{ \Gamma_{\dghA}(n), \Gamma_{\dghB}
(n) \}$ for every $n$ (as usual, we assume here that the generating
set of $G$ is fixed).
We also know that the intersection of two subgroups of finite 
index can be at most of the product of the indices. But in
fact $\Gamma_{\dghAB}$ behaves in this manner all along the way.
\begin{prop}
If $H_1,H_2 < G$ then $\Gamma_{\dghAB}(n) \leq
\Gamma_{\dghA}(n) \Gamma_{\dghB} (n)$ for every $n$.
\label{prGroups10}
\end{prop}
{\em Proof}. Let $T_1, T_2$  and $T$ be minimal right Schreier 
transversals for $H_1, H_2$ and
$H_1 \cap H_2$ respectively. Each element $g \in T$ can be
represented by a pair ($g_1, g_2$), $g_i \in T_i$, where $g_i$
is the representative of the coset $H_i g$ of $H_i$. By the
minimality of the lengths of the
coset representatives, $l(g_i) \leq l(g)$ for
each $i$. The result then follows since the pairs are 
distinct.
\hfill $\Box$ \\ \\
The proposition gives a sufficient condition for the 
intersection of two subgroups to be non-trivial.
$\Gamma_{\dghA}(n) \Gamma_{\dghB} (n) \prec  \Gamma_{\dg}(n)
\ \Lrw \
\Gamma_{\dg}(n) \not \preceq \Gamma_{\dghA}(n) \Gamma_{\dghB} (n)
\ \Lrw \ $
$\exists n \ [ \
\Gamma_{\dghA}(n) \Gamma_{\dghB} (n) < \Gamma_{\dg}(n) \ ]
\ \Lrw \
H_1 \cap H_2$ is non-trivial. For example, if $\Gamma_{\dg}(n)
\sim d^{n^e}$, where $d>1, 0 < e \leq 1$, and 
$\Gamma_{\dghA}(n), \Gamma_{\dghB} (n) \prec  \Gamma_{\dg}(n)$
then $H_1 \cap H_2$ is non-trivial.

If $H_2 < H_1 < G$ then an element $g \in G$ of minimal length 
in $H_1 g$ is also of minimal length in the coset $H_2 g$, thus
a minimal Schreier transversal for $H_1$ can be
chosen to be a subset of a minimal transversal for $H_2$.
\\
Let $1 \lrw H \stackrel{i}{\lrw} F \stackrel{\pi}{\lrw} G 
\lrw 1$, where $F$ is free and finitely generated on $X$.
Then there is a 1-1 correspondence between a
minimal Schreier transversal $T$ for $i(H)$ in $F$ and the
set of elements of the group $G$ with generating set $\{ \pi(x)
\ \mid \ x \in X \}$, given by $g=x_{i_1}\cdots x_{i_n} 
\longleftrightarrow \pi(g)=\pi(x_{i_1}) \cdots \pi(x_{i_n})$, where
$g$ and $\pi(g)$ here are already in reduced form
(note that we do not require the extension to split).
Then a (minimal) Schreier transversal for a subgroup $G'$ of
$G$ can be represented (with the above correspondence) by a 
(minimal) Schreier transversal for $\pi^{-1}(G')$ in $F$
which can
be taken to be a subset of $T$. This means that the cogrowth 
functions of subgroups of finitely generated groups
are all cogrowth functions of subgroups of finitely generated
free groups (and as will be seen in Theorem~\ref{thGroups10},
it suffices to consider the free group on 2 generators for the
equivalence classes of the cogrowth functions).
\subsection{Essential Subgroups and Their Cogrowth Functions}
We come now to essential subgroups. We call
a subgroup $H < G$ {\em essential} if
it intersects non-trivially every non-trivial subgroup of
$G$. Clearly the family $\cal E$ of essential subgroups
of a given group is a filter: if $H_1 \in \cal{E}$ and
$H_1 < H_2$ then $H_2 \in \cal{E}$, and if $H_1, H_2 \in \cal{E}$
then $H_1 \cap H_2 \in \cal{E}$.  Also a conjugate of an
essential subgroup is essential. For example, if $G$ is 
torsion-free then clearly every subgroup of finite index is
essential.
In finitely generated free groups, a subgroup is of finite
index if and only if it is finitely generated and essential.

In contrast to the situation in algebras (as will be
shown in the next section), one cannot tell whether a subgroup
of infinite index of a free group is essential or not just by
knowing the cogrowth of the subgroup.  This is due to the
``one-dimensionality'' of a subgroup generated by a single
element.
For example, let $G$ be the free group on the two generators
$x,y$ and let $H$ be the normal closure of the subgroup of $G$ 
generated by $x$. Then a minimal Schreier transversal for $H$ 
consists of all powers of $y$, hence
$\Gamma_{\dgh}(n) \sim n$. But
$H$ is not essential since $H \cap <y> = 1$. 
Thus, when $H$ is of infinite index it can be of minimal
cogrowth and still lack essentiality. On the other hand, $H$
can be essential although it has exponential cogrowth. If
$H$ is a normal subgroup of a torsion-free group $G$ then 
by definition $H$ is
essential if and only if $G/H$ is periodic (torsion). 
The well known examples of essential subgroups $H \lhd G$ of 
exponential cogrowth are when $G$ is free of rank $m
\geq 2$ and $G/H$ is the Burnside group $B(m,n)$ with $n 
\geq 665$ and odd, as shown by Adyan (\cite{Ady}).
Essential normal subgroups of intermediate growth
were constructed by Grigorchuk.

Let $\cal{CG}$ be the class of functions $\alpha (n)$ of
the following type.
$\alpha : \mbox{\bf N} \cup \{0\} \lrw \mbox{\bf N}$
is the sequence of partial sums $\sum_{i=0}^{n}f_i$ 
of the series $\sum_{i=0}^{\infty}f_i$,
such that
(i) the $f_i$-s are zero on $(r, \infty)$, where $0 \leq r$
and can be $\infty$, and positive integers otherwise,
with $f_0 = 1$;
(ii) there exists $0 < d$ such that $f_{i+1} \leq d f_i$
for every $i$.
Clearly $\cal{CG}$
includes the cogrowth functions of subgroups, but
as seen from (\ref{eqGroups10}), the set of equivalence classes
of the growth functions of groups is properly contained in
the set of equivalence classes of the members of $\cal{CG}$.
\begin{thm}
Let $G$ be the free group on $X=\{x_1, x_2\}$. Then for every
$\alpha(n) \in \cal{CG}$
there exists an essential subgroup $H$ of $G$ such that
$\Gamma_{\dgh}(n) \sim \alpha(n)$.
\label{thGroups10}
\end{thm}
{\em Proof}.
If $\alpha(n)$ is eventually-constant then any subgroup of 
finite index can be taken. So let us
assume that $\alpha(n)$ is not bounded.
We will construct in an inductive way a right
Schreier transversal $T$ for $H$. $T$ will contain two
types of sections constructed alternately:
those responsible for the desired growth (g-sections),
and those to ensure the essentiality (e-sections).
The growth function of $T$, $\Gamma_{\dt}(n)$,
will be $\preceq \Gamma_{\dgh}(n)$, which is the growth 
function of a
minimal Schreier transversal tree, but by an appropriate
definition of the coset function, and by letting the
g-sections be of sufficient depth (length) compared to the
e-sections, $T$ can be constructed such that 
$\Gamma_{\dt}(n) \sim \Gamma_{\dgh}(n)$ (if each 
g-section will be of depth equal to that of the next
e-section, the growth of $T$ will be at least as
half the growth of $\Gamma_{\dgh}$).

We start with a g-section. Here we prevent the 
occurrence of the same generator (or the same inverse of a generator)
in adjacent edges. Hence, each non-root vertex will have
either 1 or 2 out-going edges. As we will see later, an
essential section can be constructed to be of
growth $2n$, so no problem will be to achieve growth equivalent
to the minimal possible growth of $\alpha(n)$.
As for the maximal growth, $\alpha(n)$ can
grow at most as $d^n$, for some positive $d$.
So if $c$ is such that $2^c \geq d$ then 
$2^{cn} \geq d^n$, and we can construct $T$ such that
$\alpha(n) \leq \Gamma_{\dt}(cn)$. Then this 
inequality can surely be
reached when $\alpha(n)$ grows slower than $d^n$, and
thus $\Gamma_{\dt}(n)$ can be bounded by
\begin{equation}
\alpha(c^{-1}n) \leq \Gamma_{\dt}(n) \leq \alpha(2n).
\label{eqGroups33}
\end{equation}
(In general, $\Gamma_{\dt}(n)$ can be constructed so that it grows
at the fastest possible rate, but not exceeding $\alpha(n)$,
resulting in some averaging of
$\alpha(n)$, making it ``smoother'' in places of great jumps.)

Suppose we constructed the tree $T$ up to depth $p_1$, being
the first g-section. 
We label also each edge of $T$ by some $x \in X \cup X^{-1}$, 
such that the labels on the vertices are the elements one
obtains by reading off the edge labels in a path that starts
at the root and terminates at the given vertex.
We then partially define the coset function $\pi$.
If an edge labelled by $x \in X \cup X^{-1}$ goes from
the vertex $g$ to the vertex $h$ then $\pi(gx)=h$ and
$\pi(hx^{-1})=g$, that is $\pi(g)=g$ for every $g \in T$. 
Otherwise, we define $\pi$
on the set 
\begin{equation}
\{ \ gx \ \mid \ g \in T, \ l(g) < p_1 -1 \  \mbox{and} \
x \in X \cup X^{-1} \ \}
\end{equation}
in the following way. We go as far as possible on a path in the
{\em opposite} direction. That is, we start at $g$, and go
from a vertex $h'$ to a vertex
$h''$ whenever $\pi(h'x^{-1})$ is defined
and $h''=\pi(h'x^{-1})$. If the vertex $h$ is the endpoint of
this path (in general, $h$ can be $g$ itself) 
then we define $\pi(gx)=h$.
We notice that by the restriction of not labeling two adjacent
edges with the same letter, such a path as described above
will be of length $1$ (except near the root where it can
be of length $2$). The same process of defining $\pi$ will take
place at the following stages, but then, due to the e-sections,
paths as above could be of greater length. The definition of
$\pi$ on other elements of $G$ is then according to the inductive
rule $\pi(gx)= \pi(\pi(g)x)$.

Next we construct an e-section.
Assume an ordered list of the elements of $G$ is given.
We take the first element $g$ of this list,
and ``travel'' along $T$ as long as possible
with powers of $g$, starting from the root $1$ and
using the function $\pi$. If we happen to get
back to the root after some $k$-th power of $g$ then we
are done: $g^k \in H$, and we can take the next element in 
the list. Otherwise, we
will extend $T$ and $\pi$ so that $<g>$ will
intersect $H$ non-trivially. (One may look at $T$
as representing
an automaton, which has to be extended so that $g$ will be accepted
by it. Here the states of the automaton are the vertices,
with $1$ the accepting state,
the input alphabet is $X \cup X^{-1}$ and
the function is $\pi$.)
The idea is to form two paths starting from the root, one of 
a positive power of $g$ and the other of a negative power, and
to ``tie'' these paths using the function $\pi$.
Assume that when written in reduced form we have
$g = h_1 h_2 h_1^{-1}$,
where $h_2$ is of minimal length.
We start a travel from the root of $T$, this time with
$h_1$ followed by powers of $h_2$. It may still happen that we
will return to $\pi(h_1)$ after some powers of $h_2$, and in 
this case too we are done.
(We note that it is
impossible to return to the same vertex after some powers of $h_2$
before first visiting $\pi(h_1)$. This is
because by the very definition of $\pi$,
if $\pi(h_1h_2^r) = \pi(h_1h_2^s)$ then $\pi(h_1h_2^{s-r})= \pi(h_1)$.)
Otherwise, we reach 
a vertex $h \in T$ on which $\pi(hx)$ is not yet defined, where
$x$ is the next ``input'' letter of $h_1h_2^k$
for some $k \geq 0$.
We then increase $T$ by adding a path, starting with the vertex $hx$,
according to the rest of $h_1h_2^k$.
The same process of increasing the tree is then done with a path
that goes from the root with $h_1$ followed 
by powers of $h_2^{-1}$.
Since no prefix of $h_2^{-1}$ equals a suffix of $h_2$, the
two added paths must become separated one from another,
and after that happens we need not increase the tree anymore.
The endpoints of the two added paths are vertices $u_1$ and $u_2$ 
such that
$u_1 = \pi(h_1h_2^r)$ and $u_2 = \pi(h_1h_2^{-s})$,
for some $r,s \geq 0$. Assume now that
$h_2 = x w$,
where $x \in X \cup X^{-1}$, $w \in G$ and $h_2$ is written in reduced
form. Then we further extend the tree by adding a $w^{-1}$-segment 
at the vertex $u_2$ and define
$\pi(u_1 x) = u_2 w^{-1}$,
and necessarily
$\pi(u_2 w^{-1} x^{-1}) = u_1$.
Later, this construction will result (after a sufficient extension 
of $\pi$) in
\begin{equation}
\pi(g^{r+s+1})=1,
\end{equation}
i.e. $<g>$ will intersect $H$ non-trivially. Let us call $T_1$
the current tree we have.

Next we construct a g-section, extending $T_1$ and the coset
function $\pi$, according to the growth function $\alpha(n)$,
the same as before. 
If needed, we first widen the tree in the part of the
e-section to reach the desired growth, but we do not change 
$\pi$ where it is already defined.
Again, we do not define $\pi$ on the boundary of the current
tree, to ensure further increasing of the tree.
Following this stage comes an e-section with the next element
in the list. The result is a tree $T_2$. We continue in this way
indefinitely  and define $T = \bigcup_i T_i$. 
Clearly
\begin{equation}
\pi(\pi(gx)x^{-1}) = g
\end{equation}
for every $g \in T$, $x \in X \cup X^{-1}$. This makes
$T$ a right Schreier transversal for a unique subgroup $H$ of $G$,
for which $\pi$ is the function giving the coset representatives
(see \cite{Hal1}, \cite{Hal2}). $H$ is freely generated by
the non-trivial elements of the form
$gx(\pi(gx))^{-1}$,
where $g \in T$, $x \in X$ (Nielsen-Schreier theorem). 

As said before, $\Gamma_{\dt}(n)$ can be made equivalent to 
$\alpha(n)$. On the other hand, if $g$ is a
vertex in a g-section (and assume it is not on the boundary
of the section) then $l(\pi(gx)) \leq l(g)+1$, for any 
$x \in X \cup X^{-1}$. If $g$ is in an e-section then 
$l(\pi(gx)) \leq l(g)+r$, where $r$ is the length of the 
e-section. Therefore, if each g-section has at least the depth of
the next e-section, we get for any $g \in G$ 
\begin{equation}
l(\pi(g)) \leq 2l(g).
\end{equation}
Thus, the length of each element in $T$ is at most 
twice the length of the minimal element in its coset, and so
$\Gamma_{\dgh}(n) \leq \Gamma_{\dt}(2n)$ and they are equivalent
(because $\Gamma_{\dt}(n) \leq \Gamma_{\dgh}(n)$).
Combining it with (\ref{eqGroups33}) gives
\begin{equation}
\Gamma_{\dgh}(n) \sim \Gamma_{\dt}(n) \sim \alpha(n).
\end{equation}

Finally, the e-sections
make sure that for every element of
$G$ some positive power of it lies in $H$, that is $H$ is essential.
\hfill $\Box$ \\

We remark that with a little more 
effort (by adding segments corresponding
to the different group elements), the subgroups constructed in the
theorem above can have the additional property of not containing
any subgroup which is normal in $G$.

As we have seen, even when the group is torsion-free a
normal subgroup of infinite index
need not be essential although it can be of minimal cogrowth.
However, the following simple observation expresses the 
``largeness'' of normal subgroups in special cases.
If $G$ is a torsion-free
group which does not contain a non-cyclic abelian
subgroup then every two non-trivial normal subgroups of $G$ 
have non-trivial intersection.
To see it, let $1 \neq x \in H_1$, $1 \neq y \in H_2$, where 
$H_1$ and $H_2$ are normal subgroups of $G$.
Then $xyx^{-1}y^{-1} = x(yx^{-1}y^{-1}) \in H_1$ and also
$xyx^{-1}y^{-1} = (xyx^{-1})y^{-1} \in H_2$
and the result follows.
\section{Growth, Cogrowth and Essentiality in Algebras} 
Let $R$ be an associative algebra with a unit generated 
on a finite set $X$ over a field $K$.
Having the length function on the free semigroup $X^*$
generated by $X$ and the grading of $R$ by the subspaces $R^{(n)}=
\sum_{i=0}^n KX^i$, the growth function on $R$ is defined by
\begin{equation}
\Gamma_{\dr}(n) := dim \ R^{(n)}
\end{equation}
(see \cite{Row}, and also \cite{Alj} for a generalization).
The length of an element $r \in R$ is the smallest $n$ such 
that $r \in R^{(n)}$. As in groups, the equivalence class
of a growth function does not
depend upon the set of generators.

If $I$ is a right ideal of $R$ then its growth function is
$\Gamma_{\di}(n) := dim \ (I \cap R^{(n)})$,
and its cogrowth
$\Gamma_{\dri}(n) := \Gamma_{\dr}(n) - \Gamma_{\di}(n)$.
The definition can then be extended to subspaces of $R$ by
$\Gamma_{\scriptscriptstyle V}(n) := dim \ (V \cap R^{(n)})$,
and if $V \supseteq U$ then
$\Gamma_{\scriptscriptstyle{V/U}}(n) := 
\Gamma_{\scriptscriptstyle V}(n) -
\Gamma_{\scriptscriptstyle U}(n)$.

Let $V$ be a complementary subspace to $I$, that is $R=I+V$
and $I \cap V=0$. Then a basis $T$ for $V$ 
can consist of a set of (monic) monomials, which moreover is
initially  closed, that is forms a tree. Introducing a
ShortLex order on $X^*$ and extending it in the usual manner
to a partial order on $R$, such a basis can be formed from all
monomials which are minimal in their cosets (see \cite{Lew}),
resulting in
a unique minimal Schreier transversal $T$ (similar to the group
case). Then we get
\begin{eqnarray}
\Gamma_{\di}(n) &=& \mbox{card} \{ g \ \mid \  g \ \mbox{is a leading
monomial of some} \ r \in I, \ l(r) \leq n \}, \\
\Gamma_{\dri}(n) &=& \mbox{card} \{ g \in T \ \mid \ l(g) \leq n \}.
\end{eqnarray}
In fact, when $R$ is a group algebra $KG$ and $H$ a subgroup
of $G$ then the Schreier transversals for the right ideal $I$ of $R$
generated by the elements $h-1$, where $h \in H$, coincide with
the Schreier transversals for $H$ in $G$, and thus $\Gamma_{\dgh}(n)
= \Gamma_{\dri}(n)$ (keeping the generating set for $G$ fixed).

The situation concerning intersection of right ideals is
simpler than intersection of subgroups because in each $R^{(n)}$
we have intersection of finite dimensional subspaces.

If $I$ is a right ideal and $r \in R$ then
$(I : r) := \{ s \in R \ \mid \ rs \in I \}$.
\begin{prop}
If $\Gamma_{\di}(n) \not \preceq \Gamma_{\dri}(n)$ then
for every $0 \neq r \in R$, $(I : r) \neq 0$.
\label{prA77}
\end{prop}
{\em Proof}.
We may assume that $l(r) \geq 1$.
Since $\Gamma_{\di}(n) \not \preceq \Gamma_{\dri}(n)$ then
\begin{equation}
\forall C \ \exists n \ [ \ \Gamma_{\dr}(n) \geq \Gamma_{\di}(n) > 
\Gamma_{\dri}(Cn) \geq \Gamma_{\drii}(Cn) \ ].
\label{eqA10}
\end{equation}
Hence $\Gamma_{\drii}(n) \prec \Gamma_{\dr}(n)$.
Taking $C = 2l(r)$ we get $Cn \geq n+l(r)$
and by (\ref{eqA10})
\begin{equation}
\exists n \ [ \ \Gamma_{\dr}(n) > \Gamma_{\drii}(n+l(r)) \ ].
\label{eqA20}
\end{equation}
We look now at 
$\{ rg \; \mid \; g \in G, \ l(g) \leq n_0 \}$,
where $n_0$ is such that the inequality in (\ref{eqA20}) holds. 
We have here $\Gamma_{\dr}(n_0)$ elements of length
$\leq n_0+l(r)$. By (\ref{eqA20}) these elements are linearly 
dependent modulo $I$, hence there exists some 
$0 \neq s = \sum_{l(g) \leq n_0} a_g g$, $a_g \in K$, $g \in G$, 
such that $rs \in I$, i.e. $(I:r) \neq 0$.
\hfill $\Box$ \\

We now come to essentiality of right ideals. Essential right 
ideals are more common than essential subgroups.
>From a geometrical point of view, the difference is that
a right ideal generated by a single element grows in a cone-like
manner, whereas cyclic subgroups are ``1-dimensional''.
\begin{corol}
Let $R$ be without zero divisors and let $I$ be a
right ideal of $R$. If
$\Gamma_{\dri}(n) \prec \Gamma_{\dr}(n)$ then $I$ is essential.
\label{corA10}
\end{corol}
{\em Proof}.
Since $I$ is not empty it contains a right regular element.
Therefore its growth function is equivalent to the growth
function of $R$. Thus $\Gamma_{\dri}(n) \prec \Gamma_{\di}(n)$,
and by Proposition~\ref{prA77} $I$ is essential.
\hfill $\Box$  \\ \\
Remark: The converse of the above does not hold. \\

When $R$ is the group algebra $KG$, where $G$ is free of rank
$2$ and $K$ is a field, then for every $\alpha(n) \in 
\cal{CG}$ (as defined in the previous section)
there exists an essential right ideal $I$ of $R$ with
$\Gamma_{\dri}(n) \sim \alpha(n)$. For the functions which are
$\prec 2^n$ we can take the right ideal generated by the ``right
augmentation ideal'' of the subgroups constructed in 
Theorem~\ref{thGroups10}. For exponential growth we can
take, for example, the fractal ideals defined in 
\cite{Ros}.
\small


\begin{thebibliography}{99}
\bibitem{Ady} Adyan, S.I. (1975). {\em The Burnside problem
and identities in groups}. Nauka, Moscow, 1975 (Russian)
[(1979) Ergebnisse der Mathematik und ihrer Grenzgebiete, 
{\bf 95}. Springer-Verlag].
\bibitem{Alj} Aljadeff, E., Rosset, S. (1988). {\em Growth 
and uniqueness of rank}. Israel J. Math., Vol. {\bf 64}, No. 
{\bf 2}, 251-256.
\bibitem{Fab} Fabrykowski, J., Gupta, N. (1985). {\em On groups
with sub-exponential growth functions}. J. Indian Math.
Soc., Vol {\bf 49}, 249-256.
\bibitem{Gri2} Grigorchuk, R.I. (1978). {\em Symmetrical random walks
on discrete groups}. In (Ed. Dobrushin, R.L., Sinai, Ya.G.): 
Multicomponent random systems. Nauka, Moscow, 132-152
[English transl. (1980) {\em Advances in probability and related
topics}, Vol. {\bf 6}, 285-325. Marcel Dekker].
\bibitem{Gri} Grigorchuk, R.I. (1983). {\em On Milnor's problem
of group growth}. Dokl. Ak. Nauk SSSR, {\bf 271}, 31-33 (Russian)
[English transl. (1983): Soviet Math. Dokl., {\bf 28}, 23-26].
\bibitem{Gri1} Grigorchuk, R.I. (1990). {\em On growth in group
theory}. Proc. of the International Congress of
Mathematicians, Kyoto, 1990, 325-338.
\bibitem{Gro} Gromov, M. (1981). {\em Groups of polynomial growth
and expanding maps}. Publ. Math. IHES, {\bf 53}, 53-73. 
\bibitem{Hal1} Hall, M., Rad\'{o}, T. (1948). {\em On Schreier systems
in free groups}. Trans. AMS, {\bf 64}, 386-408.
\bibitem{Hal2} Hall, M. (1949). {\em Coset representations in
free groups}. Trans. AMS, {\bf 67}, 421-432.
\bibitem{Lew} Lewin, J. (1969). {\em Free modules over free algebras
and free group algebras : the Schreier technique}. Trans. AMS,
{\bf 145}, 455-465.
\bibitem{Lyn} Lyndon, R.C., Schupp, P.E. (1977).
{\em Combinatorial group theory}. Springer-Verlag.
\bibitem{Mil} Milnor, J. (1968). {\em A note on curvature and
fundamental group}. J. Diff. Geom., {\bf 2}, 1-7.
\bibitem{Mil1} Milnor, J. (1968). {\em Growth of finitely generated
solvable groups}. J. Diff. Geom., {\bf 2}, 447-449. 
\bibitem{Ros} Rosenmann, A. {\em Essentiality of 
fractal Ideals}. Inter. J. Algebra and Computation (to appear).
\bibitem{Row} Rowen, L.H. (1988) {\em Ring theory}.
Academic Press.
\bibitem{Wolf} Wolf, J.A. (1968). {\em Growth of finitely generated
solvable groups and curvature of Riemannian manifolds}.
J. Diff. Geom., {\bf 2}, 421-446. 
\end{thebibliography}
\end{document}